\documentclass[12pt]{article}
\usepackage{amsmath}
\usepackage[dvips]{graphics}

\newcommand{\no}{\nonumber}

\begin{document}

\vskip 2em \noindent \textbf{{\large Cure-rate estimation under
Case-1 interval censoring}}

\vskip 1em \noindent {\it Running title}: Cure-rate under interval
censoring

\vskip 2em \noindent {\bf Arusharka Sen}, {\it Department of
Mathematics and Statistics, Concordia University, 1455 de
Maisonneuve Boulevard West, Montreal, H3G 1M8, Canada}

\vskip 1em \noindent {\bf Fang Tan}, {\it Department of
Mathematics and Statistics, Concordia University, 1455 de
Maisonneuve Boulevard West, Montreal, H3G 1M8, Canada}

\vskip 1em \noindent \textbf{ABSTRACT. We consider nonparametric
estimation of cure-rate based on mixture model under Case-1
interval censoring. We show that the nonparametric
maximum-likelihood estimator (NPMLE) of cure-rate is non-unique as
well as inconsistent, and propose two estimators based on the
NPMLE of the distribution function under this censoring model. We
present a cross-validation method for choosing a `cut-off' point
needed for the estimators. The limiting distributions of the
latter are obtained using extreme-value theory. Graphical
illustration of the procedures based on simulated data are
provided.}


\vskip 1cm \noindent \textbf{1. Introduction}

\vskip 1em Consider a sample of individuals on each of whom some
sort of \emph{time-to-event} data is being collected, for
instance, onset time of a disease following exposure to infection,
time to death under a terminal disease, time (for criminals) to
re-offend after at least one offence etc. In most such cases,
there may be a possibility that the individual may be
\emph{immune} (e.g., not catch a disease) or get \emph{cured}
(e.g., cured of a disease or not re-offend). This is all the more
relevant when the data is subject to some kind of `open-ended'
censoring such as random censoring, double censoring or interval
censoring, where an individual being censored (i.e., event not
occurred), especially after a large amount of time, points to the
possibility of cure. In the literature, the term \emph{long-term
survival} has also been used for cure.

Cure is usually quantified by the probability of cure, or the
\emph{cure-rate}: $p=P\{X=\infty\},$ where $X$ is the
time-to-event of interest. Most of the statistical literature on
cure is based on one of the two following models for the
`improper' random variable $X$: the \emph{mixture model}, in which
$P\{X>t\}=p+(1-p)S_0(t),$ $S_0(\cdot)$ being a proper survival
function representing the finite part of $X$ (Berkson and Gage,
1952); and the \emph{bounded cumulative hazard (BCH) model}, in
which $P\{X>t\}=\exp(-\int_0^th(s)ds),$ with
$\theta:=\int_0^\infty h(s)ds<\infty,$ so that $p=\exp(-\theta)$
(see Tsodikov \emph{et al} (2003) for an excellent review).
Inference, with or without (random) censoring, has been based
mostly on either the Bayesian approach (see Yin and Ibrahim (2005)
and the references therein) or a semi-parametric approach (see
Zhao and Zhou (2006) and the references therein).

From the non-parametric point of view, it is clear that the two
models above are equivalent. Notable among the nonparametric
approaches are: Laska and Meisner (1992), who consider the NPMLE
of $p$ under random censoring when a number $m\geq 1$ of cures are
\emph{known}; Maller and Zhou (1996), who consider the value of
the Kaplan-Meier distribution function at the largest datum as an
estimator of $(1-p)$ (as is well-known, the value is less than
unity if the largest datum is censored --- an indication of cure).
See Section 2 for more comments on these two works. Another
interesting paper is Betensky and Schoenfeld (2001), who consider
a \emph{time-to-cure}, rather than just possibility of cure,
competing with time-to-event/censoring.

In this paper we study estimation of cure-rate under \emph{Case-1
interval censoring}, or \emph{current-status} data, using the
mixture model. We have been able to trace only one paper so far
under this set-up, namely Lam and Xue (2005), who work with a
semi-parametric model, allowing the cure-rate to depend on
covariates via a logit function. We consider only the parameters
$(F,p),$ the time-to-event distribution function and the
cure-rate, respectively. Of course, this is a semi-parametric
model too, but one without covariates. We show that the
Maller-Zhou idea does not work here and propose two estimators of
$p$ based on the usual (i.e., when $p=0$) NPMLE of $F$, as given
by Groeneboom and Wellner (1992). The asymptotics of the
estimators are obtained using extreme-value theory.

In Section 2, we describe the Case-1 interval censoring model with
cure-rate and show that the NPMLE of $p$ is non-unique and
inconsistent. We then propose the two estimators that depend on a
`cut-off' point. Section 3 shows how to make an optimal choice of
this cut-off point, because it involves a variance-bias trade-off
as in extremal index estimation (see, for instance, Embrechts
\emph{et al.} (1997)). In Section 4, limiting distributions of the
estimators are derived. Use of the latter to construct confidence
intervals for $p$ is straightforward.

\vskip 3em \noindent \textbf{2. Model, preliminary results and
estimators}

\vskip 1em Consider a variable of interest $X,$ say $X=$ time to
development of cancer following exposure to radiation and an
observation time $Y,$ say $Y=\mbox{ time of check-up}.$ Under
Case-1 interval censoring model, one observes the so-called
`current status' data
$$(\delta_i,Y_i), \ i=1,2,\ldots,n,\mbox{ where
}\delta_i=I(X_i\leq Y_i),$$ and $Y_1,\ldots,Y_n$ are iid with
distribution $G,$ independent of $X_1,\ldots,X_n$ which are iid
with distribution $F.$ Suppose we want to estimate $F(x)=P\{X\leq
x\}.$ The nonparametric maximum likelihood estimator (NPMLE) is
obtained by solving:
\begin{eqnarray}
& & \max_FL(F_1,\ldots,F_n)\nonumber \\ & & \mbox{subject to
}0\leq F_1\leq \ldots \leq F_n\leq 1,
\end{eqnarray}
where
$$L(F_1,\ldots,F_n)=\sum_{i=1}^n(\delta_{[i]}\log(F_i)+(1-\delta_{[i]})\log(1-F_i)),$$
and $F_i=F(Y_{(i)}),$ $Y_{(i)}:$ order-statistics for
$(Y_1,\ldots,Y_n),$ $\delta_{[i]}=$ concomitant of $Y_{(i)},$
$1\leq i\leq n.$

Solution is given by the `\emph{max-min}' formula of Groeneboom
and Wellner (1992), namely, \begin{equation} \hat{F}_i=\max_{h\leq
i}\min_{k\geq
i}\frac{\sum_{j=h}^k\delta_{[j]}}{k-h+1}.\end{equation} \noindent
\textbf{Cure-rate.} Consider again $X=$ time to cancer, this time
with \emph{possibility of no cancer} $\equiv$ \emph{cure}. Then
$X$ can be modelled as an `extended' real-valued random variable
with a \emph{defective} distribution, i.e.,
$$P(X=\infty)=p=\mbox{ \emph{cure-rate} }>0$$ so that $P(X\leq
t)=F_p(t)=(1-p)F(t)$ and $P(X>
t)=S_p(t)=p+(1-p)(1-F(t))=p+(1-p)S(t).$ In this case the
likelihood function in Eq.(1) has to be modified as $\max
L^c(p,F_1,\ldots,F_n)$ where
\begin{eqnarray} & & L^c(p,F_1,\ldots,F_n)\nonumber \\
 & = &
 \sum_{i=1}^n[\delta_{[i]}\log((1-p)F_i)+(1-\delta_{[i]})\log(p+(1-p)(1-F_i))]\nonumber\\
 & & \mbox{\emph{subject to} }0\leq p\leq 1, \ 0\leq F_1\leq \ldots \leq F_n\leq
1\nonumber\\
 & = & \sum_{i=1}^n[\delta_{[i]}\log(F_i)+(1-\delta_{[i]})\log(1-F_i)]\nonumber\\
 & & \mbox{\emph{subject to} }0\leq p\leq 1, \ 0\leq F_1\leq \ldots \leq F_n\leq
(1-p),
\end{eqnarray} writing $F_i$ for $(1-p)F_i$ in the last equality.

\vskip 1em \noindent \textbf{Failure of NPMLE.} We state the
following theorem whose proof is omitted because it is long and
technical:

\vskip 1em \noindent {\sc Theorem 1.} Let $L^c(p)=\max_{0\leq
F_1\leq \ldots \leq F_n\leq (1-p)}L^c(p,F_1,\ldots,F_n).$ Then
$$L^c(p)=L(\hat{F}_1\wedge(1-p),\ldots,\hat{F}_n\wedge(1-p)),$$
where $\wedge$ denotes `minimum' and $\hat{F}_i, \ 1\leq i\leq n,$
are as in Eq.(2).

This leads to the following two observations about the NPMLE of
$p$:

\vskip 1em \noindent {\sc Remark 1: Non-uniqueness of NPMLE.}
Obviously, $L^c(p)$ is non-increasing in $0\leq p\leq 1,$ and
$$\sup_{0\leq p\leq 1}L^c(p)=L(\hat{F}_1,\ldots,\hat{F}_n)=L^c(\hat{p}),$$
for any $0\leq \hat{p}\leq (1-\hat{F}_n).$ Hence $\hat{p}$ is
unique if and only if $(1-\hat{F}_n)=0=\hat{p}.$ This was also
observed, in the case of random censoring, by Laska and Meisner
(1992), who showed that NPMLE was unique and positive if, however,
some number $m\geq 1$ of cases of cure were known. We shall
explore this situation in a future paper.

\vskip 1em \noindent {\sc Remark 2: Non-consistency of NPMLE.}
Note that by Eq.(2),
$$\hat{F}_n=\max_{i\leq
n}\frac{\sum_{j=i}^n\delta_{[j]}}{n-i+1},$$ so that $\hat{F}_n=1$
if and only if $\delta_{[n]}=1.$ Thus for $0<p<1$ and any
$0<\varepsilon<p,$ $$P\{|\hat{F}_n-(1-p)|>\varepsilon\}\geq
P\{\hat{F}_n=1\}=P\{\delta_{[n]}=1\}=(1-p)E(F(Y_{(n)}))\rightarrow(1-p)F(\tau_G),$$
where $\tau_G=\sup\{y | G(y)=1\}.$ Hence $\hat{F}_n$ is not a
consistent estimator of $(1-p).$ This is in stark contrast to the
case of random censoring where the former was shown to be in fact
$\sqrt{n}$-consistent (asymptotically normal) by Maller and Zhou
(1996).

\vskip 1em \noindent \textbf{The proposed estimators.} Let us look
at
$$\hat{F}_n=\max_{i\leq n}\frac{\sum_{j=i}^n\delta_{[j]}}{n-i+1}=
\max_{x\leq Y_{(n)}}\frac{\sum_{j=1}^n\delta_jI(Y_j\geq
x)}{\sum_{j=1}^nI(Y_j\geq x)}.$$ Thus $\hat{F}_n$ is the maximum
of the \emph{tail-averages} of the concomitants, $\delta_{[i]}, \
1\leq i\leq n.$ Hence consider the ratio empirical process
$$p_{1n}(x):=\frac{\sum_{j=1}^n\delta_jI(Y_j\geq
x)}{\sum_{j=1}^nI(Y_j\geq x)} \rightarrow
p_1(x):=(1-p)\frac{\int_x^\infty FdG}{\int_x^\infty dG}$$ almost
surely for each $x\geq 0$ as $n\rightarrow\infty.$ Moreover, note
that
$$p_1(x)\uparrow (1-p)\mbox{ as }x\uparrow\infty$$ and
$$p_{2n}(x):=\max_{y\leq x}p_{1n}(y)\rightarrow
(1-p)\max_{y\leq x}\frac{\int_y^\infty FdG}{\int_y^\infty
dG}\uparrow (1-p)\mbox{ as }x\uparrow\infty$$ These observations
lead us to the following:

\noindent \textsc{Estimator-1.} Define
$$\hat{p}_{1n}=p_{1n}(x_n)=\frac{\sum_{j=1}^n\delta_jI(Y_j\geq
x_n)}{\sum_{j=1}^nI(Y_j\geq x_n)},$$ i.e., tail-average at a
suitable sequence $x_n\uparrow\infty$ of `cut-off' points.

Figure 1 gives a sample-plot of $p_{1n}(i)\equiv
p_{1n}(Y_{(i)})=\sum_{j=i}^n\delta_{[j]}/(n-i+1)$ against $1\leq
i\leq n,$ for $p=0.3, \ n=100.$ It is seen that for $i\approx 55,$
$p_{1n}(i)\approx 0.7=(1-p).$ For comparison, a sample-plot for
another sample with $p=0$ (i.e., \emph{no} cure) is also given.

\vskip 1em \noindent \textsc{Estimator-2.} Define
$$\hat{p}_{2n}=p_{2n}(x_n)=\max_{y\leq x_n}p_{1n}(y),$$ i.e.,
partial maximum of the tail-averages (rather than the global
maximum $\hat{F}_n$ which is inconsistent).

Figure 2 gives a sample-plot of $p_{2n}(i)=\max_{k\leq
i}\sum_{j=k}^n\delta_{[j]}/(n-i+1)$ against $1\leq i\leq n,$ for
the same sample as in Figure 1. $p_{2n}(\cdot)$ looks more stable
than $p_{1n}(\cdot),$ as is to be expected.

The choice of $x_n$ for a given sample of size $n$ is discussed in
the next section.
\newpage
\begin{figure}[!h]
\begin{center}
\scalebox{1.0}{\includegraphics{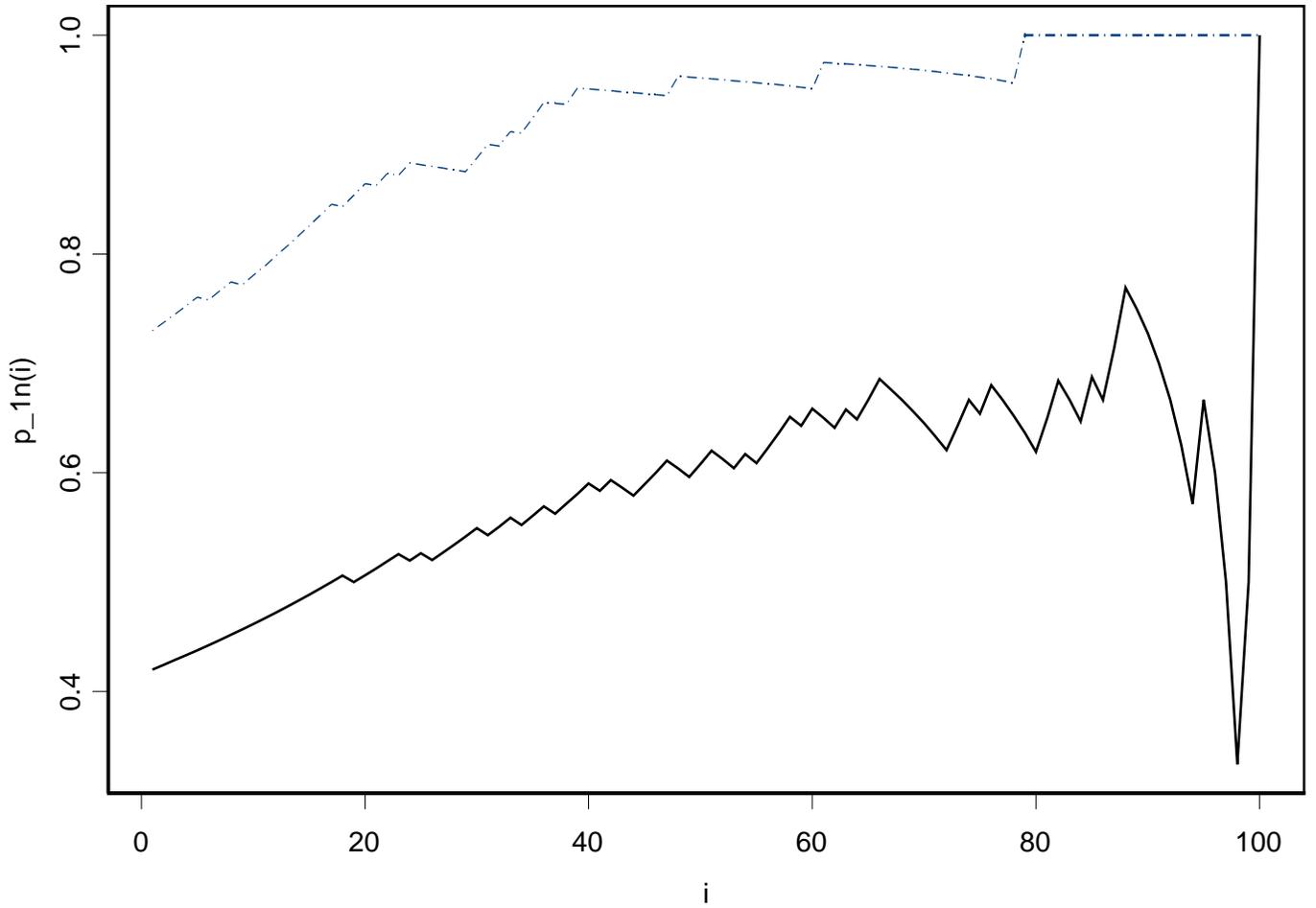}}
\end{center}
\caption{sample plot of $p_{1n}(i)$ vs. $i$: $F=$ Exp (2), $G=$
Exp (1), $n=100,$ and $p=0.3$ (\emph{solid} line), $p=0$
(\emph{broken} line).}
\end{figure}
\newpage
\begin{figure}[!h]
\begin{center}
\scalebox{1.0}{\includegraphics{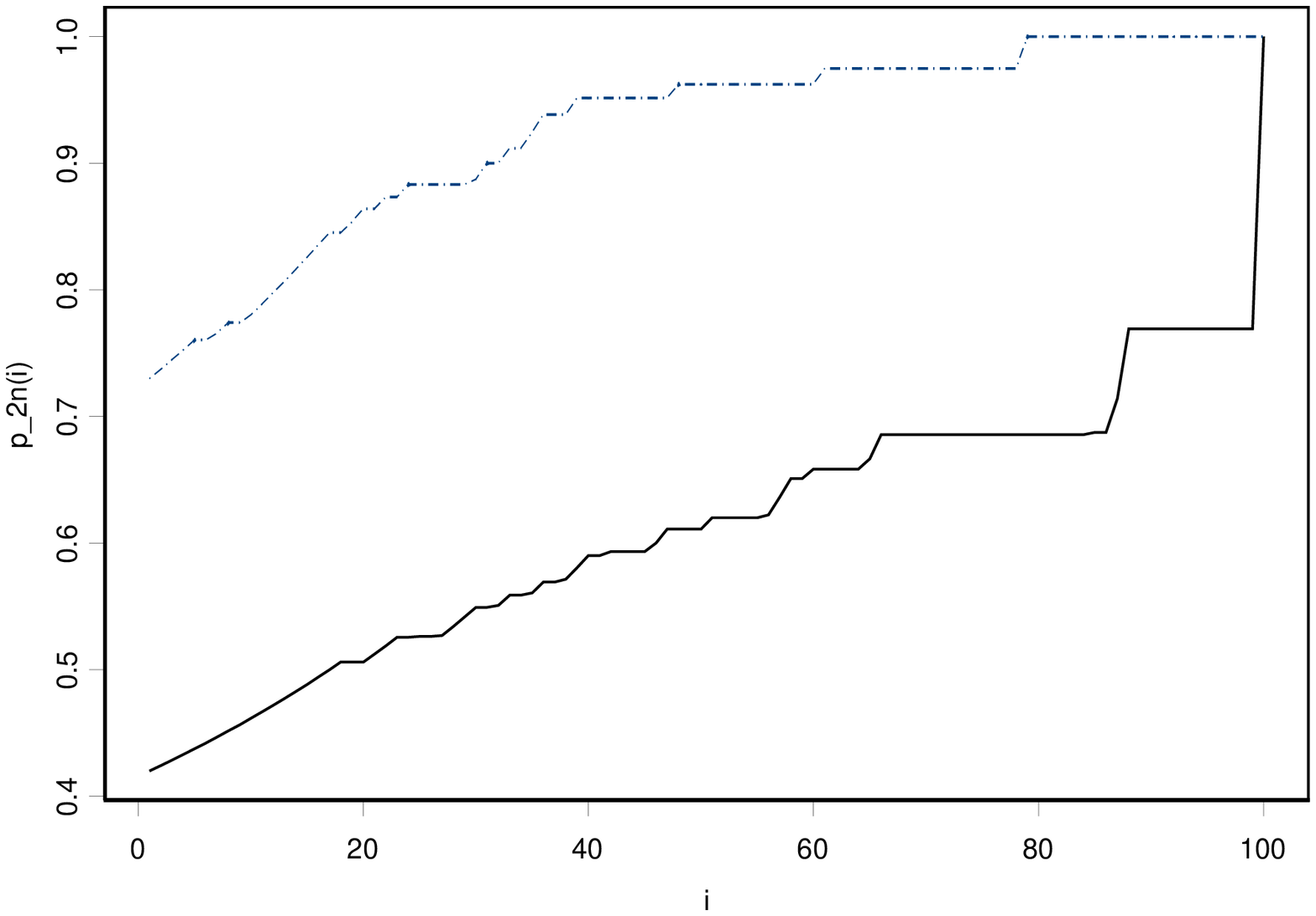}}
\end{center}
\caption{sample plot of $p_{2n}(i)$ vs. $i$: $F=$ Exp (2), $G=$
Exp (1), $n=100,$ and $p=0.3$ (\emph{solid} line), $p=0$
(\emph{broken} line).}
\end{figure}

\newpage \noindent \textbf{3. Choice of cut-off point}

Consider
\begin{eqnarray} & & \hat{p}_{1n}-(1-p)\no\\ & = &
(\hat{p}_{1n}-p_1(x_n))+(p_1(x_n)-(1-p))\no \\
& = & \frac{\sum_{j=1}^n[\delta_jI(Y_j\geq
x_n)-(1-p)(\int_{x_n}^\infty FdG/\bar{G}(x_n))I(Y_j\geq
x_n)]}{n\bar{G}(x_n)}\frac{\bar{G}(x_n)}{n^{-1}\sum_{j=1}^nI(Y_j\geq
x_n)}\no\\ & & -(1-p)\int_{x_n}^\infty (1-F)dG/\bar{G}(x_n)\no\\
 & = & A_n(x_n)C_n(x_n)- B_n(x_n),\mbox{ say,}
\end{eqnarray} where $\bar{G}(x)=\int_x^\infty dG=1-G(x).$ Now
\begin{equation}
n\bar{G}(x_n)(\mbox{var }A_n(x_n))=(1-p)\int_{x_n}^\infty
FdG/\bar{G}(x_n)-[(1-p)(\int_{x_n}^\infty
FdG/\bar{G}(x_n))]^2\rightarrow p(1-p)\end{equation} as
$x_n\rightarrow\infty.$

Further, $C_n(x_n)=O_P(1)$ (see Shorack and Wellner (1986), p.415)
and $B_n(x_n)=o(1)$ as $x_n\rightarrow\infty.$ Hence from Eq.(4),
$$\hat{p}_{1n}-(1-p)=
(\hat{p}_{1n}-p_1(x_n))+(p_1(x_n)-(1-p))=O_P((n\bar{G}(x_n))^{-1/2})+o(1),$$
as $x_n\rightarrow\infty.$

\noindent \textbf{Variance-bias trade-off.} Thus we have the
following \emph{trade-off}: as $n\rightarrow\infty,$ we must have
$x_n \uparrow \infty$ (so that the bias $-B_n(x_n)\rightarrow0$
and also $\bar{G}(x_n)\rightarrow0$), but slowly enough so that
$n\bar{G}(x_n)\rightarrow\infty$ (i.e., var
$(A_n(x_n))\rightarrow0$). A similar phenomenon occurs in the case
of the Hill estimator of extremal index in extreme value theory
(see Embrechts et al, 1997, p.341).

In view of Eq.(4)--(5), optimal order of $x_n\uparrow\infty$ could
be determined by minimizing, with respect to $x,$ the function
$$M_n(x)=(p(1-p)/n\bar{G}(x))+
(1-p)^2(\int_x^\infty(1-F)dG/\bar{G}(x))^2.$$

\vskip 1em \noindent {\sc Example 1.} Let $F, \ G$ be Exponential
($\lambda$) and Exponential ($\mu$) distributions, respectively,
i.e., $\bar{F}(x)=1-F(x)=\exp(-\lambda x),$
$\bar{G}(x)=1-G(x)=\exp(-\mu x).$ Then we have
\begin{eqnarray*} & & M_n(x)\\
&=&(p(1-p)/n\bar{G}(x))+(1-p)^2(\int_x^\infty(1-F)dG/\bar{G}(x))^2\\
&=&n^{-1}p(1-p)\exp(\mu
x)+((1-p)\mu/(\lambda+\mu))^2\exp(-2\lambda x),
\end{eqnarray*} and $(d/dx)(M_n(x))=0$ gives
$$n^{-1}p(1-p)\mu\exp(\mu x)=((1-p)\mu/(\lambda+\mu))^22\lambda\exp(-2\lambda
x),$$ or
$$x_n=(\mu+2\lambda)^{-1}\log\left(((1-p)\mu/2p\lambda(\lambda+\mu)^2)n\right).$$
Thus $n\bar{G}(x_n)=c(p,\lambda,\mu)n^{2\lambda/(\mu+2\lambda)},$
which shows that the optimal rate of convergence,\\
$(n\bar{G}(x_n))^{1/2}=O(n^{\lambda/(\mu+2\lambda)}),$ is much
slower than $\sqrt{n}.$

\vskip 1em \noindent \textbf{ Cross-validation.} Eq.(4)--(5) also
suggest that we could make a data-driven choice of $x_n$, say
$\hat{x}_n,$ as the minimizer of
$$\hat{M}_n(x):= \widehat{\mbox{ var }}(A_n(x))+\hat{B}_n^2(x)$$ with respect to
$x,$ where $\widehat{\mbox{ var }}(A_n(x))$ and $\hat{B}_n(x)$
denote suitable estimators of $\mbox{ var }(A_n(x))$ and $B_n(x),$
respectively.

Now an obvious choice of $\widehat{\mbox{ var }}(A_n(x))$ is
\begin{equation}\widehat{\mbox{ var
}}(A_n(x))=\frac{p_{2n}(x)(1-p_{2n}(x))}{\sum_{j=1}^nI(Y_j\geq
x)},\end{equation} where we have used $p_{2n}(\cdot)$ in view of
its stability, as is evident from Figure-2. The choice of
$\hat{B}_n(x),$ however, is not clear in general. Let us therefore
consider the special case of the \emph{Koziol--Green} model of
censoring:

\vskip 1em \noindent {\sc Assumption A.1.}
$1-F(x)=(1-G(x))^\alpha$ for some $\alpha>0.$

Under A.1, we have
\begin{eqnarray} B_n(x)&=&-(1-p)(1-G(x))^\alpha/(\alpha+1)\\
E(1-\delta)&=&p+(1-p)/(\alpha+1)\nonumber\\
\mbox{whence }(1-p)/(\alpha+1)&=&E(1-\delta)-p\\
\mbox{and }\alpha &=&E(\delta)/[E(1-\delta)-p].
\end{eqnarray} We then replace $E(\delta)$ by
$\bar{\delta}_n:=n^{-1}\sum_{i=1}^n\delta_i$ and $(1-p)$ by
\begin{equation}\bar{p}_{2n}:=n^{-1}\sum_{i=1}^np_{2n}(Y_i)=\int p_{2n}(x)dG_n(x),
\end{equation}
where $G_n(\cdot)$ is the empirical distribution function of
$Y_1,\ldots,Y_n.$ This is motivated as follows: for $y\geq 0,$
$$\frac{\int_y^\infty p_{2n}(x)dG_n(x)}{\bar{G}_n(y)}\approx(1-p)\frac{\int_y^\infty
(\int_x^\infty FdG/\bar{G}(x))dG(x)}{\bar{G}(y)}
=(1-p)[1-(\alpha+1)^{-2}(1-G(y))^\alpha],$$ which has bias of a
smaller order than $p_{2n}(y);$ to a first approximation, we let
$y=0$ to get $\bar{p}_{2n}.$

Thus by Eq.(6)--(10), we arrive at the following
\emph{cross-validation} function:
\begin{equation} \hat{M}^1_n(x)=\frac{p_{1n}(x)(1-p_{1n}(x))}{\sum_{j=1}^nI(Y_j\geq
x)}+(\bar{p}_{2n}-\bar{\delta}_n)^2\left[n^{-1}\sum_{j=1}^nI(Y_j\geq
x)\right]^{2\hat{\alpha}},\end{equation} where
$\hat{\alpha}=\bar{\delta}_n/(\bar{p}_{2n}-\bar{\delta}_n),$ which
could be minimized with respect to $x$ to obtain $\hat{x}_n.$

In general, motivated by Eq.(10) we could estimate the bias,
$B_n(x)=(1-p)\int_x^\infty FdG/\bar{G}(x)-(1-p),$ by
$\hat{B}_n(x):=p_{2n}(x)-\bar{p}_{2n}.$ This leads to another
cross-validation function \begin{equation}
\hat{M}^2_n(x)=\frac{p_{1n}(x)(1-p_{1n}(x))}{\sum_{j=1}^nI(Y_j\geq
x)}+(p_{2n}(x)-\bar{p}_{2n})^2 \end{equation}

Figure 3 gives sample-plots of
$\hat{M}^l_n(i)\equiv\hat{M}^l_n(Y_{(i)}), \ l=1,2.$ Both the
curves exhibit clear convex shapes with unique minima. However,
$\hat{M}^1_n(\cdot)$ shows a spurious minimum at the upper
extreme, which must be discarded. Further, the respective
minimizers are seen to underestimate $(1-p),$ so there appears to
be scope for improvement here.
\newpage
\begin{figure}[!h]
\begin{center}
\scalebox{1.0}{\includegraphics{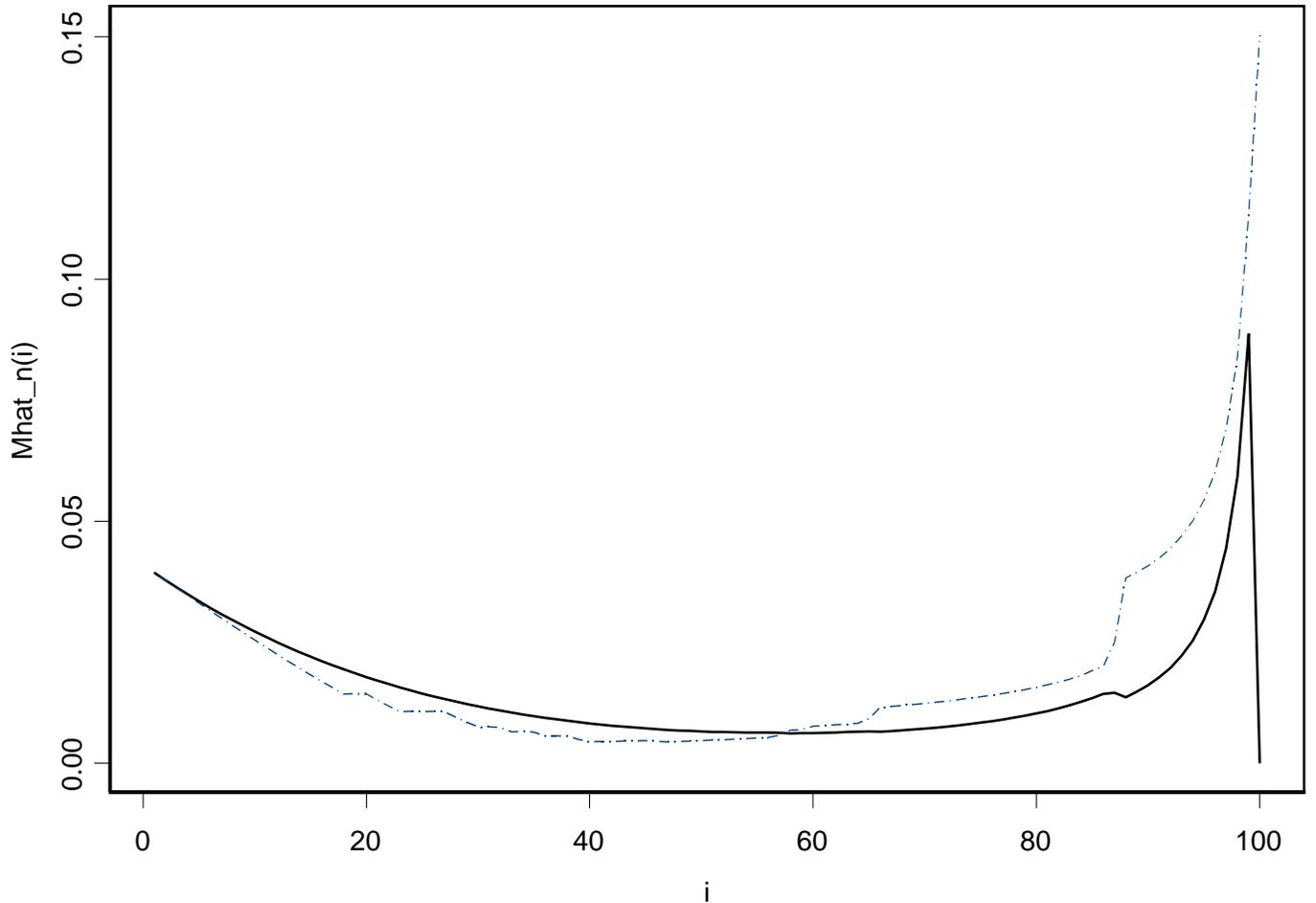}}
\end{center}
\caption{sample-plot of $\hat{M}_{n}(i)$ vs. $i$: $p=0.3,$ $F=$
Exp (2), $G=$ Exp (1), $n=100$, $\hat{M}^1_n(\cdot)$ (\emph{solid
line}: minimizer $i_1=58,$ ignoring $i=100,$ $p_{2n}(58)=0.651$),
$\hat{M}^2_n(\cdot)$ (\emph{broken line}: minimizer $i_2=47,$
$p_{2n}(47)=0.611$)}
\end{figure}

\vskip 3em \noindent \textbf{4. Limiting distributions.}

\vskip 1em Eq.(5) suggests that $\hat{p}_{1n}$ would require a
random norming, namely $(\sum_{j=1}^nI(Y_j\geq x_n))^{1/2},$ for
asymptotic normality. We establish this, as well as the limiting
distribution of $\hat{p}_{2n},$ using the asymptotic theory of
sample extremes. To this end, assume

\vskip 1em \noindent {\sc Assumption A.2.} $G(\cdot)$ belongs to
the \emph{maximum domain of attraction} of an extreme-value
distribution $G_e(\cdot),$ i.e., there exist sequences of
constants $a_n>0, \ b_n, \ n\geq 1,$ such that $G^n(a_nx+b_n)
\rightarrow G_e(x),$ or equivalently
$n\bar{G}(a_nx+b_n)\rightarrow -\log (G_e(x)),$ as
$n\rightarrow\infty,$ for each $x\in I\!\!R.$

It is well-known that, under A.2, $\sum_{j=1}^nI(Y_j\geq
a_nx+b_n)$ converges weakly to a non-homogeneous Poisson process
with mean-function $\Lambda(x)=-\log G_e(x).$ It turns out that
$\sum_{j=1}^n\delta_jI(Y_j\geq a_nx+b_n)$ converges to an
(independently) \emph{thinned} version of this process.

\vskip 1em \noindent \textsc{Lemma 1.} With
$\stackrel{d}{\rightarrow}$ denoting weak convergence in the space
$D(I\!\!R)$ of right-continuous functions on $I\!\!R$ with
left-limits, we have, as $n\rightarrow\infty,$
\begin{itemize}
\item[(a)] $N_{n}(x):=\sum_{j=1}^nI(Y_j\geq a_nx+b_n)
\stackrel{d}{\rightarrow}N(x)\equiv N([x,\infty)),$ a Poisson
process with mean $\Lambda(x)=-\log G_e(x);$ \item[(b)]
$(N_{1n}(x),N_{0n}(x))\stackrel{d}{\rightarrow} (N_1(x),N_0(x)),$
where $N_{1n}(x):=\sum_{j=1}^n\delta_jI(Y_j\geq a_nx+b_n),$
$N_{0n}(x):=\sum_{j=1}^n(1-\delta_j)I(Y_j\geq a_nx+b_n),$ and
$N_1(x), N_0(x)$ are \emph{independent} Poisson processes with
mean-functions $\mu_1(x):=(1-p)\Lambda(x), \
\mu_0(x):=p\Lambda(x)$ respectively. \item[(c)] Further,
$N_1(x)\stackrel{d}{=}\sum_{j=1}^{N(x)}\eta_j,$ and
$N_0(x)\stackrel{d}{=}\sum_{j=1}^{N(x)}(1-\eta_j),$where
$(\eta_1,\eta_2,\ldots)\mbox{ are iid Bernoulli }(1-p),$
independent of $N(\cdot),$ and $N(\cdot)$ is the Poisson process
defined in Part (a) above.
\end{itemize}
\noindent {\sc Proof:}

\noindent (a) This is a classical result. For a proof see, for
instance, Embrechts et al. (1997).

\noindent (b) First, consider weak convergence of $N_{1n}(x)$
alone. It is enough to verify convergence of the
finite-dimensional distributions
$(N_{1n}(x_1),\ldots,N_{1n}(x_k)), \ k\geq 1$ (see, for instance,
Karr (1991), Theorem 1.21, p.14). For the sake of convenience let
us consider just two points, $(x_1,x_2)$ with $x_1<x_2.$ Then with
$i=\sqrt{-1}$ and any real numbers $t_1,t_2,$
\begin{eqnarray*} & & E[\exp(it_1N_{1n}(x_1)+it_2N_{1n}(x_2))]\\
&=& \left(E[\exp(\delta_1\{it_1I(Y_1\geq a_nx_1+b_n)+it_2I(Y_1\geq
a_nx_2+b_n)\})]\right)^n\\
&=&
\left[(p+(1-p)\int_0^\infty(1-F)dG+(1-p)\int_0^{a_nx_1+b_n}FdG)
+e^{it_1}(1-p)\int_{a_nx_1+b_n}^{a_nx_2+b_n}FdG \right.\\& &
\left.+ e^{it_1+it_2} (1-p)\int_{a_nx_2+b_n}^\infty FdG\right]^n\\
&=&\left[1+n^{-1}n\bar{G}(a_nx_2+b_n)\left\{(1-p)(e^{it_1}-1)
\int_{a_nx_1+b_n}^{a_nx_2+b_n}FdG/\bar{G}(a_nx_2+b_n)\right.\right.\\
& & \left.\left.+(1-p)(e^{it_1+it_2}-1) \int_{a_nx_2+b_n}^\infty
FdG/\bar{G}(a_nx_2+b_n)\right\}\right]^n\\
&\rightarrow&
\exp\left((1-p)\Lambda(x_2)\left\{(e^{it_1}-1)(\Lambda(x_1)\Lambda^{-1}(x_2)-1)
+(e^{it_1+it_2}-1)\right\}\right),
\end{eqnarray*} whence the result. Note that here we have used the fact
that as $n\rightarrow\infty,$ $(a_nx+b_n)\rightarrow\tau_G,$ so
that $\int_{a_nx+b_n}^\infty FdG/\bar{G}(a_nx+b_n)\rightarrow 1.$
The joint weak convergence of $(N_{1n}(x),N_{0n}(x)),$ as well as
their asymptotic independence, follow by exactly similar
arguments.

\noindent (c) The representations of $(N_{1}(x),N_{0}(x))$ are
obvious.$\fbox{}$

Next note that
\begin{equation} p_{1n}(x_n)=\sum_{j=1}^n\delta_jI(Y_j\geq
x_n)/\sum_{j=1}^nI(Y_j\geq
x_n)=N_{1n}(x'_n)/N_n(x'_n),\end{equation} where
$x'_n=(x_n-b_n)/a_n.$ Therefore, in addition to the weak
convergence in Lemma 1, we need \emph{strong approximation} by a
Poisson process. This follows in a straightforward way from
Einmahl (1997) and is stated below:

\vskip 1em \noindent \textsc{Theorem 2.} Under A.2, on some
probability space one can construct the random variables
$(\delta_i,Y_i), \ i=1,2,\ldots,$ and a sequence of Poisson
processes $N'_n=(N'_{1n},N'_{0n})$ on $I\!\!R\times I\!\!R$, where
for each $n\geq 1,$ $N'_{1n},\ N'_{0n}$ are \emph{independent}
with
mean-functions $\mu_1(x),\ \mu_0(x),$ respectively, such that as $n\rightarrow\infty,$\\
$\sup_{x:0<G_e(x)<1}|N_{1n}(x)-N'_{1n}(x)|\stackrel{P}{\rightarrow}0,$\\
$\sup_{x:0<G_e(x)<1}|N_{0n}(x)-N'_{0n}(x)|\stackrel{P}{\rightarrow}0.$

\vskip 1em \noindent {\sc Proof:} Follows by arguments similar to
the proof of Corollary 2.6, p.37, of Einmahl (1997).$\fbox{}$

We are now ready to state the limiting distributions of our
estimators. In Theorem 3 below, by `$\lim$' we mean \emph{limit in
distribution}.

\vskip 1em \noindent {\sc Theorem 3.} Under A.2, if
$n\bar{G}(x_n)\rightarrow\infty$ as $n\rightarrow\infty,$ then
\begin{itemize} \item[(a)] $\Lambda(x'_n)\rightarrow\infty,$
where $x'_n=(x_n-b_n)/a_n;$ \item[(b)] let $$Z_{1n}=
\frac{(\sum_{j=1}^nI(Y_j\geq
x_n))^{1/2}(p_{1n}(x_n)-(1-p))}{\sqrt{p(1-p)}};$$ then
\begin{eqnarray*} & & \lim_{n\rightarrow\infty}Z_{1n}\\
&=&\lim_{n\rightarrow\infty}\sqrt{N(x'_n)}\left[\frac{\sum_{j=1}^{N(x'_n)}\eta_{j}}{N(x'_n)}
-(1-p)\right]/\sqrt{p(1-p)} = \mbox{ Normal }(0,1),\end{eqnarray*}
where $(\eta_1,\eta_2,\ldots)\mbox{ are iid Bernoulli }(1-p)$ as
in Lemma 1, Part (c); \item[(c)] let
$$Z_{2n}=\frac{(\sum_{j=1}^nI(Y_j\geq
x_n))^{1/2}(p_{2n}(x_n)-(1-p))}{\sqrt{p(1-p)}};$$ then
\begin{eqnarray*} & &\lim_{n\rightarrow\infty}Z_{2n}\\
&=&\lim_{n\rightarrow\infty}\sqrt{N(x'_n)}\sup_{x\leq
x'_n}\left[\frac{\sum_{j=1}^{N(x)}\eta_{j}}{N(x)}
-(1-p)\right]/\sqrt{p(1-p)} = \mbox{ half-Normal }(0,1),
\end{eqnarray*} where `half-Normal' $(0,1)$ is the distribution of
$|\mbox{ Normal }(0,1)|.$
\end{itemize}

\vskip 1em \noindent {\sc Proof:}

\noindent (a) Since extreme-value distributions are all
continuous, the convergence $|G^n(a_nx+b_n) - G_0(x)|\rightarrow
0$ is uniform in $x.$ Now $n\bar{G}(x_n)\rightarrow\infty
\Rightarrow G^n(x_n)=G^n(a_nx'_n+b_n)\rightarrow 0,$ hence
$G_0(x'_n)\rightarrow 0.$ The result follows because
$\Lambda(x'_n)=-\log G_0(x'_n).$

\noindent (b) Note that \begin{eqnarray*} & &
\lim_{n\rightarrow\infty}\frac{(\sum_{j=1}^nI(Y_j\geq
x_n))^{1/2}(p_{1n}(x_n)-(1-p))}{\sqrt{p(1-p)}}\\
&=&\lim_{n\rightarrow\infty}(N_{1n}(x'_n)+N_{0n}(x'_n))^{1/2}
\left[\frac{N_{1n}(x'_n)}{N_{1n}(x'_n)+N_{0n}(x'_n)}
-(1-p)\right]/\sqrt{p(1-p)}\\ &=&
\lim_{n\rightarrow\infty}(N'_{1n}(x'_n)+N'_{0n}(x'_n))^{1/2}
\left[\frac{N'_{1n}(x'_n)}{N'_{1n}(x'_n)+N'_{0n}(x'_n)}
-(1-p)\right]/\sqrt{p(1-p)},\end{eqnarray*} by Theorem 2. The
result now follows using the representation in Lemma 1, Part (c),
and the \emph{random central limit theorem}, since
$(\eta_1,\eta_2,\ldots)\mbox{ are iid Bernoulli }(1-p),$
independent of $N(\cdot),$ and further, by Part (a) above,
$\Lambda(x'_n)\rightarrow\infty,$
$N(x'_n)/\Lambda(x'_n)\stackrel{P}{\rightarrow}1,$ as
$n\rightarrow\infty.$

\noindent (c) This result too follows as in Part (b) above, by
noting that $$\lim_{n\rightarrow\infty}\sqrt{N(x'_n)}\sup_{x\leq
x'_n}\left[\frac{\sum_{j=1}^{N(x)}\eta_{j}}{N(x)}
-(1-p)\right]/\sqrt{p(1-p)}=\lim_{n\rightarrow\infty}\sqrt{n}\sup_{m\geq
n}\left[\frac{\sum_{j=1}^m\eta_{j}}{m}
-(1-p)\right]/\sqrt{p(1-p)}.$$ Weak convergence of the sequence on
right-hand-side to the half-Normal distribution is established in
Robbins et al (1968) (see also Stute (1983) for a generalization
to $M$-estimators).$\fbox{}$

\vskip 1em \noindent {\sc Remark 1.} Figures 4 and 5 give
histograms of $Z_{1n}$ and $Z_{2n},$ respectively, based on 5000
samples each. Either of $Z_{1n}$ and $Z_{2n}$ may easily be used
to construct confidence intervals for $(1-p).$ However, note that
limiting variance of $Z_{1n}=1>1-2\pi^{-1}=$ limiting variance of
$Z_{2n}.$ Hence the latter may be a better choice. On the other
hand, Figure-5 shows that the convergence of $Z_{2n}$ to the
half-Normal distribution is \emph{not} very good.
\newpage
\begin{figure}[!h]
\begin{center}
\scalebox{1.0}{\includegraphics{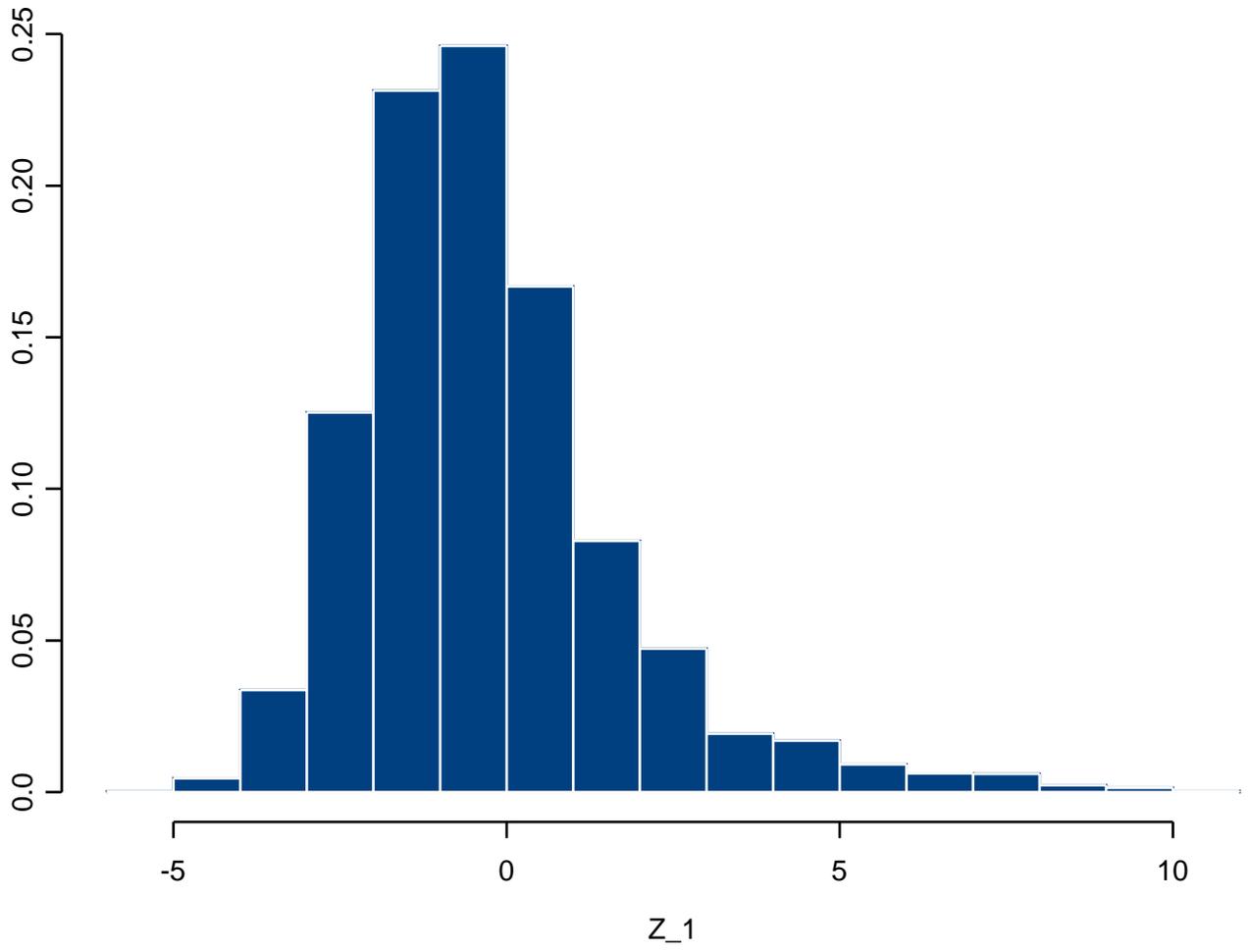}}
\end{center}
\caption{histogram of studentized $Z_{1n}$: $p=0.3,$ $F=$ Exp (2),
$G=$ Exp (1), $n=100$, based on 5000 samples}
\end{figure}

\begin{figure}[!h]
\begin{center}
\scalebox{1.0}{\includegraphics{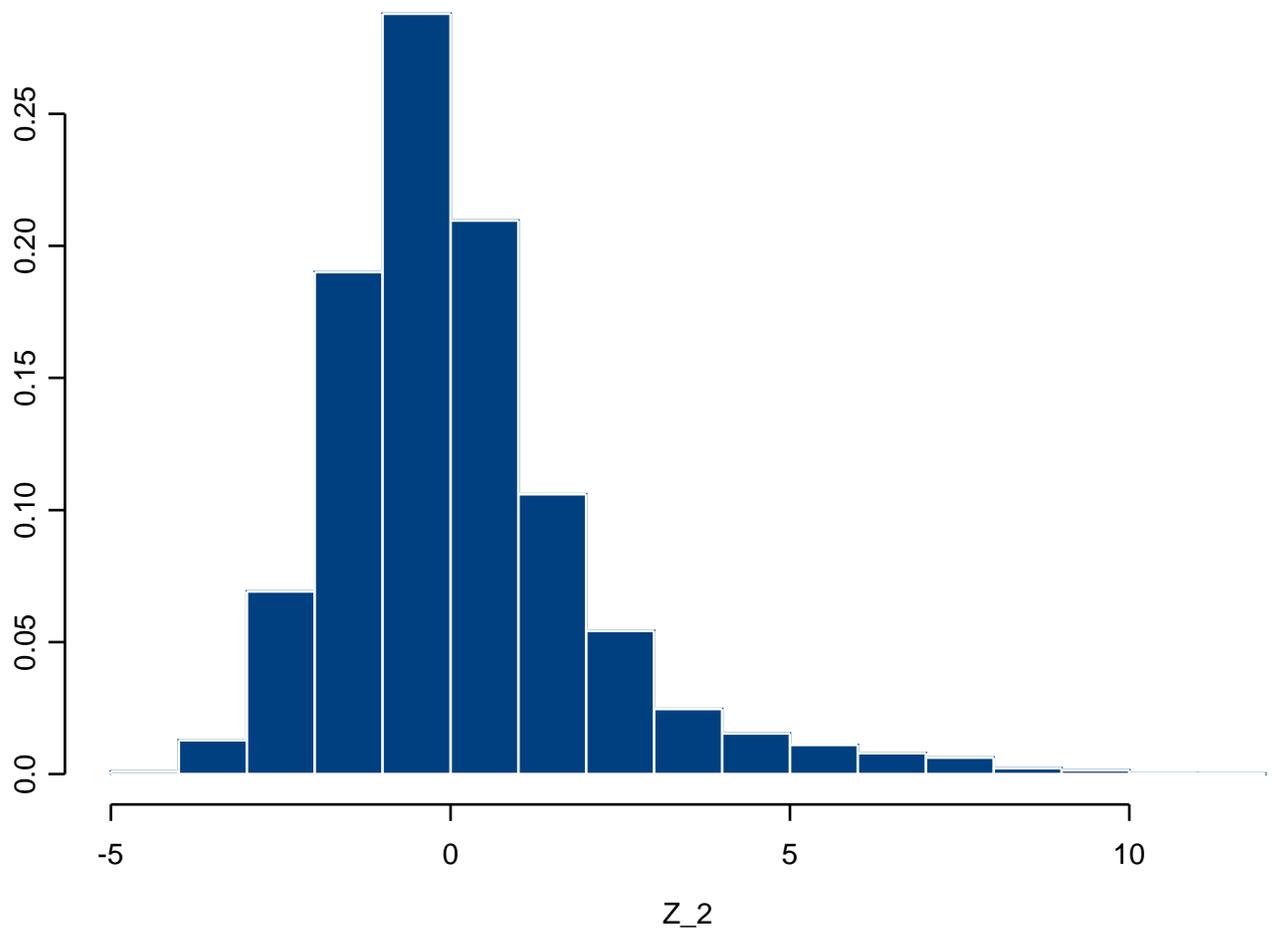}}
\end{center}
\caption{histogram of studentized $Z_{2n}$: $p=0.3,$ $F=$ Exp (2),
$G=$ Exp (1), $n=100$, based on 5000 samples}
\end{figure}
\newpage
\noindent \textit{Acknowledgements}: This research was supported
by a Discovery grant of NSERC, Canada, to A.Sen.

\vskip 3em {\centerline {\bf REFERENCES}}
\begin{description}
\item Berkson, J. \& Gage, R.P. (1952). Survival curve for cancer
patients following treatment. \emph{J. Amer. Statist. Assoc.}
\textbf{47}, 501-515. \item Betensky, R.A. \& Schoenfeld, D.A.
(2001). Nonparametric estimation in a cure model with random cure
times. \emph{Biometrics} \textbf{57}, 282--286. \item Embrechts,
P., Kl\"{u}ppelberg, C. \& Mikosch, T. (1997). \emph{Modelling
extremal events. For insurance and finance.} Springer-Verlag,
Berlin. \item Groeneboom, P. \& Wellner, J.A. (1992).
\emph{Information bounds and nonparametric maximum likelihood
estimation.} Birkh\"{a}user Verlag, Basel. \item Karr, A.F.
(1991). \emph{Point processes and their statistical inference},
2nd ed. Marcel Dekker, New York. \item Lam, K. F. \& Xue, H.
(2005). A semiparametric regression cure model with current status
data. \emph{Biometrika} \textbf{92}, 573--586. \item Laska, E.M.
\& Meisner, M.J. (1992). Nonparametric estimation and testing in a
cure rate model. \emph{Biometrics} \textbf{48}, 1223-1234. \item
Maller, R.A. \& Zhou, X. (1996). \emph{Survival analysis with
long-term survivors}. Wiley, Chichester. \item Robbins, H.,
Siegmund, D. \& Wendel, J. (1968). The limiting distribution of
the last time $s\sb{n}\geq n\varepsilon $. \emph{Proc. Nat. Acad.
Sci.} \textbf{61}, 1228--1230. \item Shorack, G.R. \& Wellner,
J.A. (1986). \emph{Empirical processes with applications to
statistics}. Wiley, New York. \item Stute, W.(1983). Last passage
times of $M$-estimators. \emph{Scand. J. Statist.} \textbf{10},
301--305. \item Tsodikov, A. D., Ibrahim, J. G. \& Yakovlev, A. Y.
(2003). Estimating cure rates from survival data: an alternative
to two-component mixture models. \emph{J. Amer. Statist. Assoc.}
\textbf{98}, 1063--1078. \item Yin, G. \& Ibrahim, J.G. (2005).
Cure rate models: a unified approach. \emph{Canad. J. Statist.}
\textbf{33}, 559--570. \item Zhao, X. \& Zhou, X. (2006).
Proportional hazards models for survival data with long-term
survivors. \emph{Statist. Probab. Lett.} \textbf{76}, 1685--1693.
\end{description}

\end{document}